\documentclass[11pt]{amsart}
\usepackage{amsmath,amsfonts,amscd,amssymb,epsf}
\newtheorem{theorem}{Theorem}[section]

\newtheorem{lemma}[theorem]{Lemma}
\newtheorem{corollary}[theorem]{Corollary}
\theoremstyle{definition}

\theoremstyle{remark} \newtheorem{remark}[theorem]{Remark}
\numberwithin{equation}{section}
\newcommand{\field}[1]{\ensuremath{\mathbb{#1}}}

\newcommand{\RR}{\field{R}}

\DeclareMathOperator{\id}{id}

 \DeclareMathOperator{\Diff}{Diff}

\DeclareMathOperator{\PSL}{PSL} \DeclareMathOperator{\PSU}{PSU}

\newcommand{\lo}{{\mathbb{L}}}
\newcommand{\Z}{{\mathbb{Z}}}
\newcommand{\U}{{\mathbb{U}}}
\newcommand{\R}{{\mathbb{R}}}
\newcommand{\Del}{{\mathbb{D}}}
\newcommand{\C}{{\mathbb{C}}}
\newcommand{\bk}{\backslash}
\newcommand{\pa}{\partial}

\newcommand{\ov}{\overline}

\newcommand{\z}{\bar{z}}

\begin{document}
\title{Bers isomorphism on the universal Teichm\"uller curve }
\date{\today}
\author{Lee-Peng Teo}\address{Faculty of Information Technology\\Multimedia Univerity\\
Jalan Multimedia\\ 63100, Cyberjaya\\Selangor, Malaysia}
\email{lpteo@mmu.edu.my} 
\begin{abstract}
We study the Bers isomorphism between the Teichm\"uller space of
the parabolic cyclic group and the universal Teichm\"uller curve.
We prove that this is a group isomorphism and its derivative map
gives a remarkable relation between Fourier coefficients of cusp
forms and Fourier coefficients of vector fields on the unit
circle. We generalize the Takhtajan-Zograf metric to the
Teichm\"uller space of the parabolic cyclic group, and prove that
up to a constant, it coincides with the pull back of the
Velling-Kirillov metric defined on the universal Teichm\"uller
curve via the Bers isomorphism.
\end{abstract}

\maketitle

\section{Introduction}

Let $\U$ be the upper half plane and $G$ a Fuchsian group. We
denote by $T(G)$ the Teichm\"uller space of $G$ and $F(G)$ the
Bers fiber space of $G$.
 In \cite{Bers}, Bers proved the following result.
\begin{theorem}\label{Bers}
Let $G$ be a torsion free Fuchsian group, $a$ a point in $\U$,
$\hat{a}$ the image of $a$ under the natural projection $\U
\rightarrow \U/G$. Let $\dot{G}$ be another torsion free Fuchsian
group and $h : \U/\dot{G} \rightarrow (\U/G) \setminus \{\hat{a}\}$
a conformal bijection. Then there is a canonical biholomorphic
bijection $T(\dot{G}) \rightarrow F(G)$.
\end{theorem}

In case $G=\{\id\}$ is the trivial group, Bers theorem implies a
biholomorphic bijection $T(\Gamma_0) \rightarrow F(1)$, where
 $\Gamma_0$ is the parabolic subgroup $\{ z \mapsto z+n\;
|\; n \in \Z \}$ of $\PSL(2, \RR)$, and $F(1)$ is the universal
Teichm\"uller curve. Both the universal Teichm\"uller space $T(1)$
and the universal Teichm\"uller curve $F(1)$ have group
structures. We are going to show in Section \ref{section2} that as
a subspace of $T(1)$, $T(\Gamma_0)$ is also a subgroup. Moreover,
we prove in Section \ref{Section3} that the Bers isomorphism on
$T(\Gamma_0)$ is a group isomorphism.

In \cite{Teo}, we proved that the universal Teichm\"uller curve
$F(1)$ is isomorphic to the homogenuous space
$S^1\bk\text{Homeo}_{qs}(S^1)$, where $\text{Homeo}_{qs}(S^1)$ is
the space of quasi-symmetric homeomorphisms of the unit circle.
Using this isomorphism, we prove in  Section \ref{section4} that
the derivative of the Bers isomorphism at the origin of
$T(\Gamma_0)$ gives a remarkable relation between the Fourier
coefficients of cusp forms on the upper half plane and Fourier
coefficients of vector fields on the unit circle.

In \cite{Teo}, we defined the Velling-Kirillov metric on the
universal Teichm\"uller curve $F(1)$, which we proved is the
unique right invariant K\"ahler metric on $F(1)$. In \cite{TZ},
Takhtajan and Zograf defined a new K\"ahler metric on $T_{g,n}$,
the Teichm\"uller space of Riemann surfaces of genus $g$ with $n$
punctures, which plays an important role in the local index
theorem. In Section \ref{section2}, we generalize the definition
of this metric to $T(\Gamma_0)$ and still call it the
Takhtajan-Zograf metric. By definition, this metric is
right-invariant. Hence it is not surprising that we have the
following theorem (see Theorem \ref{metric}).
\begin{theorem}\label{maintheorem}
The pull back of the Velling-Kirillov metric on $F(1)$ via the
biholomorphic bijection $T(\Gamma_0) \rightarrow F(1)$ is a
multiple of the Takhtajan-Zograf metric on $T(\Gamma_0)$.
\end{theorem}

\section{Teichm\"uller spaces and K\"ahler
metrics}\label{section2}
In this section, we review the results we need from Teichm\"uller
theory and our previous paper \cite{Teo}. The readers can consult
standard texts in Teichm\"uller theory \cite{Lehto, Nag} or our
papers \cite{Teo, TT2} for further details.

Let $\U$, $\lo$, $\Del$ and $\Del^*$ be respectively the upper
half plane, the lower half plane, the unit disc and the exterior
of unit disc. In this paper, we are going to use both the upper
half plane model and the unit disc model for the universal
Teichm\"uller space $T(1)$. In case of the disc model, it is
described as follows. Let $L^{\infty}(\mathbb{D})$ be the complex
Banach space of bounded Beltrami differentials on $\mathbb{D}$. We
denote by $L^{\infty}(\mathbb{D})_1$ the unit ball of
$L^{\infty}(\mathbb{D})$. For $\mu \in L^{\infty}(\mathbb{D})_1$,
we extend $\mu$ to $\mathbb{D}^*$ by reflection
\begin{align*}
\mu(z) = \ov{\mu\left(\frac{1}{\z}\right)}\frac{z^2}{\z^2},
\hspace{2cm} z\in \mathbb{D}^*.
\end{align*}
There is a unique quasi-conformal mapping $w_{\mu}$ fixing $-1,
-i, 1$ which solves the Beltrami equation
\begin{align*}
(w_{\mu})_{\z} = \mu (w_{\mu})_z.
\end{align*}
It fixes the unit circle $S^1$, $\mathbb{D}$ and $\mathbb{D}^*$.
The universal Teichm\"uller space is defined as a set of
equivalence classes of normalized quasi-conformal mappings
\begin{align*}
T(1) = L^{\infty}(\mathbb{D})_1/\sim,
\end{align*}
where $\mu \sim \nu$ if and only if $w_{\mu} = w_{\nu}$ on the
unit circle. The tangent space at identity is identified with
$\Omega^{-1,1}(\mathbb{D})$, the space of harmonic Beltrami
differentials on $\mathbb{D}$. More explicitly, let
\begin{align*}
A_{\infty}(\mathbb{D})= \left\{\phi \;\text{holomorphic on }
\mathbb{D}\; \Bigr\vert\; \sup_{z\in \mathbb{D}}
|(1-|z|^2)^2\phi(z)| < \infty \right\}
\end{align*}
be the space of bounded holomorphic quadratic differentials on
$\mathbb{D}$. Then
\begin{align*}
\Omega^{-1,1}(\mathbb{D})=\left\{\left.-\frac{(1-|z|^2)^2}{2}
\ov{\phi(z)}\; \right\vert\; \phi \in
A_{\infty}(\mathbb{D})\right\}.
\end{align*}  The upper half plane model for the
Teichm\"uller space $T(1)\simeq L^{\infty}(\U)_1/\sim$ is defined
analogously where we replace $\mathbb{D}$ and $\Del^*$ in the
definition above by $\U$ and $\lo$ respectively; the reflection on
$\mu$ is defined by $\mu(z)=\ov{\mu(\z)}$; the quasi-conformal
mapping $w_{\mu}$ is normalized to fix the points $0,1, \infty$
and two quasi-conformal mappings $w_{\mu}$ and $w_{\nu}$ are
equivalent if and only if they agree on the extended real line
$\hat{\R}=\R\cup\{\infty\}$.

The universal Teichm\"uller curve $F(1)$ is defined as
\begin{align*}
F(1)=\Bigl\{ ([\mu], z) \;\Bigr\vert \;[\mu] \in T(1), z \in
f^{\mu} (\mathbb{D})\Bigr\}.
\end{align*}
Here $f^\mu$ is the quasiconformal map that is holomorphic on
$\mathbb{D}^*$, normalized such that $f(\infty) = \infty$,
$f'(\infty) =1$, $(f-z) (\infty) = 0$ and satisfies the Beltrami
equation $f_{\z} = \mu f_z$ on $\mathbb{D}$. The tangent space at
identity is identified with $\Omega^{-1,1}(\mathbb{D}) \oplus \C$.

Let $\text{Homeo}_{qs}(S^1)$ be the space of  quasi-symmetric
homeomorphisms on the unit circle, i.e. it consists of
orientation preserving homeomorphism of the unit circle
$\omega:S^1\rightarrow S^1$ that satisfies
\begin{align}\label{def1}\frac{1}{M}\leq \frac{\omega(e^{2\pi
i(x+t)})-\omega(e^{2\pi i x})}{\omega(e^{2\pi ix})-\omega(e^{2\pi
i (x-t)})} \leq M, \hspace{1cm}\forall x,t\in \R,
0<t<\frac{1}{4},\end{align}for some constant $M$. We identify
$S^1\bk\text{Homeo}_{qs}(S^1)$ as the subspace of
$\text{Homeo}_{qs}(S^1)$ consisting of those $\omega$ that fix the
point $1$. Every $\omega\in S^1\bk\text{Homeo}_{qs}(S^1)$ can be
extended to be a symmetric (with respect to $S^1$ )
quasi-conformal mapping which we still denote by $\omega$. It has
a unique conformal welding $\omega=g^{-1}\circ f$ satisfying the
conditions 1) $f$ and $g$ are quasi-conformal mappings. 2)
$\left.f\right|_{\Del^*}$ and $\left. g\right|_{\Del}$ are
conformal and depend only on $\left.\omega \right|_{S^1}$. 3)
$f(\infty)=\infty, f'(\infty)=1, g(0)=0$. We call
$(\left.f\right|_{\Del^*},\left. g\right|_{\Del} )$ the pair of
univalent functions associated to the point $\omega\in
S^1\bk\text{Homeo}_{qs}(S^1)$, which we also write as $(f,g)$. We
identify $T(1)$ as the subspace of $S^1\bk\text{Homeo}_{qs}(S^1)$
consisting of those $\omega$ whose corresponding $f$ satisfies the
additional condition $(f-z)(\infty)(0)$, i.e. the power series
expansion of $f$ has the form \begin{align}\label{seriesf}f(z)=z +
\frac{b_1}{z}+\frac{b_2}{z^2}+\ldots,\hspace{1cm}z\in
\Del^*.\end{align} In fact, for every $\mu \in L^{\infty}(\Del)_1$
extended by symmetry to $\Del^*$, there exists a unique symmetric
quasi-conformal mapping $\omega$ that satisfies the Beltrami
equation $\omega_{\z} =\mu \omega_z$ and whose corresponding $f$
in the conformal welding of $\omega$ has the form \eqref{seriesf}.
We denote this $\omega$ by $\omega_{\mu}$ and the corresponding
conformal welding by $(f^{\mu}, g_{\mu})$. The isomorphism
$\mathfrak{T}:F(1)\xrightarrow{\sim} S^1\bk\text{Homeo}_{qs}(S^1)$
is given by
$$F(1)\ni ([\mu],z)\mapsto \left.\sigma_w\circ \omega_{\mu}
\right|_{S^1}\in S^1\bk\text{Homeo}_{qs}(S^1),$$ where
$w=g_{\mu}^{-1}(z)\in \Del$ and
\begin{align}\label{sigma}\sigma_w(\zeta)=\frac{1-\bar{w}}{1-w}\frac{\zeta-w}{1-\zeta\bar{w}}\in
\PSU(1,1).\end{align}

For the purpose of next section, we introduce the space
\begin{align*}
\mathcal{Q}=\left\{ \omega_{\mu} \;|\; \mu\in
L^{\infty}(\Del)_1\right\}.
\end{align*}
From our reasoning above, it is isomorphic to $L^{\infty}(\Del)_1$
via the canonical correspondence $\mu \mapsto \omega_{\mu}$.
Moreover, there is a real analytic isomorphism $\mathfrak{J} :
\mathcal{Q}/\sim \times \Del \rightarrow F(1)$ given by
$$\mathcal{Q}/\sim \times \Del \ni ([\omega_{\mu}], z) \mapsto ([\mu],
g_{\mu}(z))\in F(1),$$where the equivalence relation $\sim$ on
$\mathcal{Q}$ is defined in the obvious way: $\omega_{\mu}\sim
\omega_{\nu}$ if and only if $\omega_{\mu}$ and $\omega_{\nu}$
agree on $S^1$. We also define $\mathfrak{K}:\mathcal{Q}\times
\Del\rightarrow S^1\bk\text{Homeo}_{qs}(S^1)$ by $(\omega_{\mu},
z)\mapsto \left.\sigma_z\circ \omega_{\mu}\right|_{S^1}$.
Obviously, $\mathfrak{K}=\mathfrak{T}\circ \mathfrak{J}\circ \Pi$,
where $\Pi:\mathcal{Q}\times \Del\rightarrow
\mathcal{Q}/\sim\times\Del$ is the canonical projection.

 The tangent space
at identity of $S^1\bk\text{Homeo}_{qs}(S^1)$ is identified with
the Zygmund class vector fields\footnote{See \cite{Teo} for
definition.} on $S^1$. Given a one parameter flow $\omega_t\in
S^1\bk\text{Homeo}_{qs}(S^1)$ where $\omega_0 = \id$, it defines
the tangent vector $$\upsilon =u(e^{i\theta})\frac{\pa}{\pa
\theta} = \sum_{n \in \Z} c_n e^{i n \theta}\frac{\pa}{\pa
\theta}$$where $u$ is defined by $\dot{\omega}(z) = \frac{d}{dt}
\Bigr\vert_{t=0} \omega_t (z) = iz u(z)$. Since
$\omega_t(S^1)=S^1$ and $\omega_t(1)=1$ for all $t$, we find that
$$\hspace{1cm}c_{-n} = \bar{c}_{n}\hspace{1cm}\text{and}\hspace{1cm} \sum_{n\in \Z} c_n
=0.$$ The Velling-Kirillov metric is a right invariant metric on
$S^1\bk\text{Homeo}_{qs}(S^1)$. At the origin, it is given by
$$ \Vert \upsilon \Vert_{VK}^2 = \sum_{n=1}^{\infty} n|c_n|^2.
$$

The group structure on $T(1)\simeq L^{\infty}(\U)_1/\sim$ is
induced by the composition of quasi-conformal mappings, i.e.
$[\mu]*[\nu]=[\lambda]$ if and only if $w_{\mu}\circ w_{\nu} \sim
w_{\lambda}$; whereas the group structure on $F(1)$ is induced
from the group structure on $S^1\bk\text{Homeo}_{qs}(S^1)$, which
is defined by composition of quasi-symmetric homeomorphisms.

Let $\Gamma_0$ be the parabolic cyclic group $\Gamma_0=\{ z\mapsto
z+n \; \bigr\vert \; n\in \Z\}$. The Teichm\"uller space
$T(\Gamma_0)$ is defined as
$$L^{\infty}(\U, \Gamma_0)_1/\sim,$$ where $L^{\infty}(\U,
\Gamma_0)_1$ is the subspace of $L^{\infty}(\U)_1$ consisting of
those $\mu$ satisfying
$$\mu(z+1)=\mu(z), \hspace{1cm}z\in \U$$ and $\sim$ is the same
equivalence relation defined on $L^{\infty}(\U)_1$.

Given $[\mu] \in T(\Gamma_0)$, let $\Gamma_{\mu}=w_{\mu}\circ
\Gamma_0\circ w_{\mu}^{-1}$. It is a parabolic cyclic subgroup of
$\PSL(2,\R)$ generated by $$\gamma_0 = w_{\mu}\circ \theta_0\circ
w_{\mu}^{-1}, \hspace{1cm}\theta_0(z)=z+1.
$$
Putting $z=0$ and $z\rightarrow \infty$ in the relation $\gamma_0
\circ w_{\mu}= w_{\mu}\circ \theta_0,$ we find that
$\gamma_0(0)=1$ and $\gamma_0(\infty)=\infty$. Consequently,
$\gamma_0=\theta_0$, $\Gamma_{\mu}=\Gamma_0$ and
$w_{\mu}(z+1)=w_{\mu}(z)+1$. As a result, $T(\Gamma_0)$ can be
equivalently defined as equivalence classes of symmetric
quasi-conformal mappings $w_{\mu}$ satisfying the condition
$w_{\mu}(z+1)=w_{\mu}(z)+1$. Then it is easily seen that
$T(\Gamma_0)$ is a subgroup of $T(1)$.

\begin{remark}\label{isom}
It is well known that the universal Teichm\"uller space $T(1)$ is
isomorphic to the space of quasi-symmetric homeomorphisms on
$\hat{\R}$ fixing $0, 1, \infty$, where a homeomorphism $u
:\R\rightarrow \R$ is called quasi-symmetric if and only if $u$ is
monotonically increasing and there exists a constant $M$ such that
\begin{align*}
\frac{1}{M}\leq \frac{u(x+h)-u(x)}{u(x)-u(x-h)}\leq M\hspace{1cm}
\forall x,h\in \R, h>0.
\end{align*}The isomorphism is given by $[\mu]\mapsto
\left.w_{\mu}\right|_{\R}$. The subspace $T(\Gamma_0)$ is then
identified with the space of quasi-symmetric homeomorphisms on
$\R$ that fix all the integers and satisfy the condition
$u(x+1)=u(x)+1$ for all $x$.
\end{remark}

The tangent space at the origin of $T(\Gamma_0)$ is identified
with
$$\Omega^{-1,1}(\U, \Gamma_0)=\{ -2y^2 \ov{\phi}\;:\; \phi\in
A_{\infty}(\U, \Gamma_0)\},$$ where $A_{\infty} (\U ,\Gamma_0)$ is
the space of cusp forms
\begin{align*}
A_{\infty} (\U ,\Gamma_0) = \bigl\{ \phi : \U \rightarrow \C \;
\text{holomorphic}\; \bigr\vert \sup_{z \in \U} \vert y^2 \phi(z)
\vert < \infty, \phi(z+1) =\phi(z)\bigr\}.
\end{align*}

In \cite{TZ}, Takhtajan and Zograf defined a new K\"ahler metric
on $T_{g,n}$, the Teichm\"uller space of Riemann surfaces of genus
$g$ with $n$ punctures. We generalize this metric to be a right
invariant metric on $T(\Gamma_0)$. At the origin, it is defined by
\begin{align}\label{TZmetric1}
\left\langle \mu,\nu \right\rangle_{TZ} = \iint\limits_F
\mu(z)\ov{\nu(z)} d^2z, \hspace{1cm}\mu,\nu\in\Omega^{-1,1}(\U,
\Gamma_0) .
\end{align}
Here $F=\{ x+iy \;\vert \; 0<x<1, 0<y<\infty \}$ is a fundamental
domain of $\Gamma_0$ on $\U$. By definition,
$\mu=-2y^2\bar{\phi}$, where $\phi\in A_{\infty} (\U ,\Gamma_0) $
and therefore has an expansion
\[
\phi(z) = \sum_{n=1}^{\infty} \alpha_n e^{2\pi i n z}.
\]
A direct computation gives
\begin{align}\label{TZmetric2}
\left\Vert \mu \right\Vert_{TZ}^2 =\frac{3}{32 \pi^5}
\sum_{n=1}^{\infty} \frac{|\alpha_n|^2}{n^5}.
\end{align}

\section{Bers Isomorphism between $T(\Gamma_0)$ and
$F(1)$.}\label{Section3} In this section, we are going to give a
separate proof of the Bers isomorphism $\mathfrak{B}:
T(\Gamma_0)\rightarrow F(1)$. The basic idea follows the same line
as in \cite{Bers}. But we use a different model for $F(1)$.
Besides, the proof is much simpler in this case. On the other
hand, we are also going to show that $\mathfrak{B}$ is in fact a
group isomorphism.

Let $p(z)=e^{2\pi i z}$. Then $p: \U \rightarrow \Del\setminus
\{0\}$ is a holomorphic covering map which can be extended to be a
continuous map $p: \U\cup\R \rightarrow \Del\cup S^1\setminus
\{0\}$. Given $\mu \in L^{\infty}(\U, \Gamma_0)_1$, let $\nu$ be
defined such that \begin{align}\label{relation1}\nu \circ p
\frac{\ov{p'}}{p'} =\mu.\end{align} This is well defined in view
of the periodicity of $\mu$. We define the pre-Bers mapping $B:
L^{\infty}(\U, \Gamma_0)_1\rightarrow
 \mathcal{Q}\times \Del$ by $\mu \rightarrow (\omega_{\nu}, \omega_{\nu}(0))$.
 Then there exists a unique holomorphic covering map $p_{\mu}:
 \U \rightarrow \Del\setminus \{w_{\nu}(0)\}$
 such that the following diagram
 \begin{align}\label{comm1}
\begin{CD}
\U  @> w_{\mu} >> \U \\
@VV p V     @VV p_{\mu} V \\
 \mathbb{D}\setminus \{0\} @> \omega_{\nu} >> \mathbb{D}
 \setminus \{\omega_{\nu}(0) \}
\end{CD}
\end{align}
is commutative, i.e.
\begin{align}\label{commutative}
p_{\mu}\circ w_{\mu}(z) = w_{\nu}\circ p(z),\hspace{1cm} z\in \U.
\end{align}
This equation can be extended continuously to $\U\cup \R$.
 By uniqueness of holomorphic mappings, we have
 $p_{\mu}=\gamma_{\mu}\circ p$, where $\gamma_{\mu}\in \PSU(1,1)$
 is a fractional linear transformation that maps $0$ to
 $\omega_{\nu}(0)$. Moreover, putting $z=0$ in
 \eqref{commutative}, we find that $\gamma_{\mu}(1)=1$. This
 determine $\gamma_{\mu}$ uniquely as\begin{align*}
\gamma_{\mu}(\zeta)=
\sigma_{\omega_{\nu}(0)}^{-1}(\zeta)=\frac{(1-\omega_{\nu}(0))\zeta+
\omega_{\nu}(0)(1-\ov{\omega_{\nu}(0)})}{(1-\ov{\omega_{\nu}(0)})+\zeta
(1-\omega_{\nu}(0))\ov{\omega_{\nu}(0)}}
 \end{align*}
with $\sigma$ given by \eqref{sigma}. Therefore, we have
\begin{align}\label{relation2}p\circ w_{\mu}(z) = \sigma_{w_{\nu}(0)}\circ w_{\nu}\circ
p(z),\hspace{1cm}z\in \U\cup\R.\end{align}Define
$\mathfrak{P}:L^{\infty}(\U,\Gamma_0)_1\rightarrow
S^1\bk\text{Homeo}_{qs}(S^1)$ to be $\mathfrak{K}\circ B$.
Explicitly, it is given by
$$\mathfrak{P}(\mu) = \left. \sigma_{w_{\nu}(0)}\circ w_{\nu}\right|_{S^1}.$$
In other words, $\mathfrak{P}$ is defined by the commutative
diagram
\begin{align}\label{comm2}
\begin{CD}
\U \cup\R @> w_{\mu} >> \U\cup\R \\
@VV p V     @VV p V \\
 \mathbb{D}\cup S^1\setminus \{0\} @> \sigma_{w_{\nu}(0)}\circ\omega_{\nu} >>
 \mathbb{D}\cup S^1
 \setminus \{0\}
\end{CD}
\end{align}

\begin{lemma}\label{lemma1}
If $\mu_1,\mu_2\in L^{\infty}(\U,\Gamma_0)_1$ and $\mu_1\sim
\mu_2$, then $\mathfrak{P}(\mu_1)=\mathfrak{P}(\mu_2)$.
\end{lemma}
\begin{proof}
If $\mu_1\sim \mu_2$, then $\left.w_{\mu_1}\right|_{\R} =
\left.w_{\mu_2}\right|_{\R}$. From \eqref{relation2}, it is easily
deduced that then $\left. \sigma_{w_{\nu_1}(0)}\circ
w_{\nu_1}\right|_{S^1}=\left. \sigma_{w_{\nu_2}(0)}\circ
w_{\nu_2}\right|_{S^1}$, i.e.,
$\mathfrak{P}(\mu_1)=\mathfrak{P}(\mu_2)$.

\end{proof} From this lemma, we find that the map $\mathfrak{P}$ descends to the map
$\mathcal{B}: L^{\infty}(\U,\Gamma_0)_1/\sim=T(\Gamma_0)
\longrightarrow S^1\bk\text{Homeo}_{qs}(S^1)$. Finally, the Bers
isomorphism $\mathfrak{B}:T(\Gamma_0)\rightarrow F(1)$ is given by
$\mathfrak{B}= \mathfrak{T}^{-1}\circ \mathcal{B}$. Working out
explicitly,
$$T(\Gamma_0) \ni [\mu] \xrightarrow{\hspace{0.6cm}\mathfrak{B}\hspace{0.6cm}} ([\nu],
f^{\nu}(0))\in F(1).$$ From \eqref{relation1} and the holomorphic
dependance of $f^{\nu}(0)$ on $\nu$, we conclude immediately that
$\mathfrak{B}$ is a holomorphic mapping.

On the other hand, by Remark \ref{isom}, definition \eqref{def1},
the identification $F(1)\simeq S^1\bk \text{Homeo}_{qs}(S^1)$ and
the commutative diagram \eqref{comm2} restricted to maps between
$\R$ and $S^1$, we see that Bers isomorphism is nothing but a
correspondence between  quasi-symmetric homeomorphisms on the unit
circle that fix the point $1$ and the lifting of these
homeomorphisms to the universal cover of $S^1$, i.e. $\R$, under
the covering map $p: \R \rightarrow S^1$. It immediately follows
that

\begin{lemma}
The mapping $\mathcal{B} :T(\Gamma_0)\rightarrow
S^1\bk\text{Homeo}_{qs}(S^1)$ is a bijection.
\end{lemma}

Moreover, from the commutative diagram
\begin{align*}
\begin{CD}
\U \cup\R @> w_{\mu_1} >> \U\cup\R @> w_{\mu_2} >> \U\cup\R\\
@VV p V     @VV p V  @VV p V\\
 \mathbb{D}\cup S^1\setminus \{0\} @> \sigma_{w_{\nu_1}(0)}\circ\omega_{\nu_1} >>
 \mathbb{D}\cup S^1
 \setminus \{0\} @> \sigma_{w_{\nu_2}(0)}\circ\omega_{\nu_2} >>
 \mathbb{D}\cup S^1
 \setminus \{0\}
\end{CD}
\end{align*}
we conclude immediately that

\begin{lemma}
The mapping $\mathcal{B} :T(\Gamma_0)\rightarrow
S^1\bk\text{Homeo}_{qs}(S^1)$ is a group homomorphism.
\end{lemma}
Grouping together, we have shown that
\begin{theorem}
The Bers isomorphism $\mathfrak{B}: T(\Gamma_0)\rightarrow F(1)$
is a biholomorphism between complex manifolds and a group
isomorphism.
\end{theorem}

\section{The derivative mapping}\label{section4}
In this section, we study the derivative of the Bers isomorphism
at the origin. The following theorem gives a significance relation
between the Fourier coefficients of cusp forms and Fourier
coefficients of vector fields on $S^1$.

\begin{theorem}\label{isom1}
The derivative of the mapping $\mathfrak{P}: T(\Gamma_0)
\rightarrow  \text{Homeo}_{qs}(S^1)/S^1$ at the origin is given by
the following linear isomorphism
\begin{align*}
D_0 \mathfrak{P} : \Omega^{-1,1}(\U, \Gamma_0) &\longrightarrow
T_0
(S^1\bk\text{Homeo}_{qs}(S^1))\\
 -2y^2 \sum_{n=1}^{\infty}\ov{ \alpha_n  \exp(2 \pi i n z)}
&\mapsto
\frac{i}{4\pi^2}\Biggl(\sum_{n=1}^{\infty}\frac{\alpha_n}{n^3}
e^{in\theta} -\sum_{n=1}^{\infty} \frac{\bar{\alpha}_n}{n^3}
e^{-in \theta}\Biggr)+c_0.
\end{align*}
where $$c_0=\frac{1}{4\pi^2
i}\left(\sum_{n=1}^{\infty}\frac{\alpha_n}{n^3}-\sum_{n=1}^{\infty}
\frac{\bar{\alpha}_n}{n^3}\right).$$
\end{theorem}
\begin{proof}
Given $\mu \in \Omega^{-1,1}(\U, \Gamma_0)$, for $t$ in a
neighborhood of $0$, we consider the one parameter family
$w_{t\mu}$. Let $v_{t\mu}= \sigma_{\omega_{t\nu}(0)}\circ
\omega_{t\nu}$, where $\mu$ and $\nu$ are related by
\eqref{relation1} and by definition
$\mathfrak{P}([t\mu])=\left.v_{t\mu}\right|_{S^1}$. From
\eqref{comm2}, we have
\begin{align*}
 v_{t\mu}\circ p(z)=p\circ w_{t\mu}(z) , \hspace{1cm} z\in \U\cup
\R.
\end{align*}
It is well known that (see e.g. \cite{Ahlfors2, Nag}) $w_{t\mu}$
is real analytic in $\U$ and therefore $v_{t\mu}$ is real analytic
in $\Del$. Taking derivative with respect to $t$ and putting
$t=0$, we have \begin{align}\label{equation1} \dot{v}_{\mu}\circ
p(z)=p'(z) \dot{w}_{\mu}(z), \hspace{1cm} z\in \U,\end{align}
where $\dot{w}_{\mu}=\left.\tfrac{d}{dt}\right|_{t=0}w_{t\mu}$ and
$\dot{v}_{\mu}=\left.\tfrac{d}{dt}\right|_{t=0}v_{t\mu}$. By
continuity, \eqref{equation1} still holds on $\R$. Let
$$\Phi_{\mu}(z)=\frac{i}{8\pi^3}\sum_{n=1}^{\infty}\frac{\alpha_n}{n^3}\exp(2\pi
i n z)$$ so that $\mu(z) =
-2y^2\ov{\Phi_{\mu}^{\prime\prime\prime}(z)}$. By a well-known
theorem of Ahlfors \cite{Ahlfors1}, we have
\begin{align}\label{equation2}
\dot{w}_{\mu}(z) =
\frac{(z-\z)^2}{2}\ov{\Phi_{\mu}^{\prime\prime}(z)}+(z-\z)\ov{\Phi_{\mu}^{\prime}(z)}+\ov{\Phi_{\mu}(z)}+\Phi_{\mu}(z)
+p(z)+\ov{p(z)},
\end{align}
where $p(z)$ is a polynomial of degree two. Since $w_{t\mu}$ fixes
$0,1,\infty$, we have $\dot{w}_{\mu}$ vanishes at $0,1,\infty$,
and we find that $p(z)$ is a constant and is equal to $p(z)=p(0)=
-\Phi(0).$ Let \begin{align}\label{equationv}\sum_{n\in \Z} c_n
e^{in\theta} =D_0\mathfrak{P}(\mu) =
\left.\left(\frac{\dot{v}_{\mu}}{iz}\right)\right|_{S^1}.\end{align}
Restricted to $\R$, equation \eqref{equation1} and
\eqref{equation2} give us
\begin{align*}
i \exp(2\pi i x)\sum_{n\in \Z} c_n \exp(2\pi i nx)  =& 2\pi i
\exp(2\pi i x)\left(\ov{\Phi_{\mu}(x)}+\Phi_{\mu}(x)
+p(x)+\ov{p(x)}\right).
\end{align*}
Therefore,
\begin{align*}
\sum_{n\in \Z} c_n \exp(2\pi i
nx)=\frac{i}{4\pi^2}\sum_{n=1}^{\infty}
\frac{\alpha_n}{n^3}\exp(2\pi i
nx)-\frac{i}{4\pi^2}\sum_{n=1}^{\infty}
\frac{\bar{\alpha}_n}{n^3}\exp(-2\pi i nx)+c_0,
\end{align*}
where
$$c_0=-2\pi\left(\Phi(0)+\ov{\Phi(0)}\right).$$
This implies the assertion of the theorem.

\end{proof}

\begin{remark}
In \cite{Teo}, we proved that if $\upsilon =\sum_{n\in \Z} c_n
e^{in\theta}\in T_0(S^1\bk\text{Homeo}_{qs}(S^1))$, then
$$\sum_{n=1}^{\infty} n^{\alpha} |c_n|^2 <\infty$$ for all $\alpha<2.$ Theorem
\ref{isom1} then implies the series $\sum_{n=1}^{\infty}
|\alpha_n|^2n^{-s} $ converges absolutely for all
$\text{Re}\;s>4$, which is a well-known result if
$\phi=\sum_{n=1}^{\infty}\alpha_n \exp(2\pi i n z)$ is a cusp form
of a cofinite Fuchsian group.
\end{remark}

\begin{theorem}
The derivative of the Bers isomorphism $\mathfrak{B}: T(\Gamma_0)
\rightarrow F(1)$ at the origin is given by the following linear
isomorphism
\begin{align*}
D_0 \mathfrak{B} : \Omega^{-1,1}(\U, \Gamma_0) &\longrightarrow
\Omega^{-1,1}(\Del)\oplus \C\\
 -2y^2 \sum_{n=1}^{\infty}\ov{ \alpha_n  \exp(2 \pi i n z)}
&\mapsto \left(
\frac{(1-|z|^2)^2}{8\pi^2}\sum_{n=2}^{\infty}\frac{n^3-n}{n^3}\bar{\alpha}_n\z^{n-2},
-\frac{\bar{\alpha}_1}{4\pi^2}\right)
\end{align*}
\end{theorem}
\begin{proof}
Given $\mu\in \Omega^{-1,1}(\U, \Gamma_0)$, by definition,
$D_0\mathfrak{B}(\mu)=(\lambda, a)$, where $\lambda$ is the
projection of $\nu$ defined by \eqref{relation1} to harmonic
Beltrami differentials, and
$a=\dot{f}^{\nu}(0)=\left.\tfrac{d}{dt}\right|_{t=0}f^{t\nu}(0).$
We can compute $\lambda$ by using the formula of projection.
Instead, we reason as follows. Let $\omega_{t\lambda}$ be the
corresponding family of quasi-conformal mappings defined by
$t\lambda$, for $t$ in a small neighborhood of $0$. Then by
definition,
$\left.\dot{\omega}_{\lambda}\right|_{S^1}=\left.\dot{\omega}_{\nu}\right|_{S^1}$,
where as usual,
$\dot{\omega}_{\lambda}=\left.\tfrac{d}{dt}\right|_{t=0}\omega_{t\lambda}$
and
$\dot{\omega}_{\nu}=\left.\tfrac{d}{dt}\right|_{t=0}\omega_{t\nu}$.
By definition,
$\lambda=-\tfrac{(1-|z|^2)^2}{2}\ov{\phi_{\lambda}}$ for some
$\phi_{\lambda}\in A_{\infty}(\Del)$. Let $$\Phi_{\lambda}(z) =
\sum_{n=2}^{\infty} \beta_n
z^{n+1},\hspace{1cm}\phi_{\lambda}(z)=\sum_{n=2}^{\infty}(n^3-n)
\beta_n z^{n-2}$$ so that $\Phi_{\lambda}^{\prime\prime\prime}
=\phi_{\lambda}.$ By the same formula of Ahlfors \eqref{equation2}
applied to the disc model, we have
\begin{align*}
\dot{\omega}_{\lambda}(z)
=-\frac{(1-|z|^2)^2}{2}\ov{\Phi_{\lambda}^{\prime\prime}(z)}-z(1-|z|^2)\ov{\Phi_{\lambda}^{\prime}(z)}
-z^2\ov{\Phi_{\lambda}(z)}+\Phi_{\lambda}(z) +q(z)
\end{align*}
where $q(z)= a_0+a_1 z+a_2z^2$, $a_0= -\bar{a}_2$ and $a_1$ is
purely imaginary. Restricted to $S^1$, we have
\begin{align}\label{equationomega}
\dot{\omega}_{\nu}(z) = \dot{\omega}_{\lambda}(z) =
-\sum_{n=2}^{\infty}\bar{\beta}_n z^{1-n} + \sum_{n=2}^{\infty}
\beta_n z^{n+1} + a_0+a_1 z+a_2z^2.
\end{align}
 On the other hand, from the relation
$\omega_{t\nu}=g_{t\nu}^{-1}\circ f^{t\nu}$, we have
\begin{align}\label{equation3}\dot{\omega}_{\nu}(z)=
-\dot{g}_{\nu}(z)+\dot{f}^{\nu}(z), \hspace{1cm}z\in
\C.\end{align} From the fact that $g_{t\nu}(0)=0$ and the form of
$f^{t\nu}$ given by \eqref{seriesf}, we find that $a_0=0$ and
hence $a_2=0$. Equation \eqref{equation3} also gives us
$\dot{f}^{\nu}(0)= \dot{w}_{\nu}(0)$.

Now using the definition $v_{t\mu}=\sigma_{\omega_{t\nu}(0)}\circ
\omega_{t\nu}$ and \eqref{sigma}, we have
\begin{align*}
\dot{v}_{\mu}(z) =& \left(\left.
\frac{d}{dt}\right|_{t=0}\sigma_{\omega_{t\nu}(0)}\right)(z) +
\dot{\omega}_{\nu}(z)\\
=&-\dot{\omega}_{\nu}(0)+(\dot{\omega}_{\nu}(0)-\ov{\dot{\omega}_{\nu}(0)})
z+\ov{\dot{\omega}_{\nu}(0)}z^2+\dot{\omega}_{\nu}(z).
\end{align*}
Comparing this equation with equations \eqref{equationv},
\eqref{equationomega} and the result of Theorem \ref{isom1}, we
find that
\begin{align*}
\dot{w}_{\nu}(0)=&
-i\bar{c}_1=-\frac{1}{4\pi^2}\bar{\alpha}_1\hspace{0.5cm}\text{and
for $n\geq 2$,}\hspace{0.5cm} \beta_n =\; ic_n =
-\frac{1}{4\pi^2}\frac{\alpha_n}{n^3}.
\end{align*}
The assertion of the theorem follows.
\end{proof}
\begin{remark}
$T(\Gamma_0)$ can be thought of as the parameter space of all the
one-punctured surfaces. This theorem implies that the tangent
vector on $T(\Gamma_0)$ that purely 'moves the puncture'
corresponds to the automorphic form which has only the first
coefficient nonzero. On the other hand, the automorphic form with
vanishing first coefficient does not 'move the puncture'.
\end{remark}

\section{Takhtajan-Zograf metric and Velling-Kirillov metric}
In \cite{Teo} and \cite{TT2}, we have studied some properties of
the Velling-Kirillov metric on the universal Teichm\"uller curve
$F(1)$. We can pull back this metric to $T(\Gamma_0)$ via the Bers
isomorphism, which we still call the Velling-Kirillov metric. A
straight-forward computation gives
\begin{theorem}\label{metric}
On $T(\Gamma_0)$, the Velling-Kirillov metric is a multiple of the
Takhtajan-Zograf metric. More precisely, for $\mu_1, \mu_2 \in
\Omega^{-1,1}(\U,\Gamma_0)$,
\begin{align*}
\langle \mu_1, \mu_2\rangle_{VK} =
\frac{2\pi}{3}\iint\limits_{\Gamma_0\bk\U}
\mu_1(z)\ov{\mu_2(z)}d^2z =\frac{2\pi}{3}\langle \mu_1,
\mu_2\rangle_{TZ}.
\end{align*}
\end{theorem}
\begin{proof}
By definition and Theorem \ref{isom1},
\begin{align*}
\left\Vert \mu\right\Vert_{VK}^2 =\sum_{n=1}^{\infty} n\left|
\frac{i}{4\pi^2}\frac{\alpha_n}{n^3}\right|^2 =
\frac{1}{16\pi^4}\sum_{n=1}^{\infty}\frac{|\alpha_n|^2}{n^5}.
\end{align*}
Compare to the definition \eqref{TZmetric1} of the
Takhtajan-Zograf metric and \eqref{TZmetric2}, the result follows.
\end{proof}
Since we have proved in \cite{Teo} that the Velling-Kirillov
metric is K\"ahler, we obtain immediately
\begin{corollary}
The Takhtajan-Zograf metric on $T(\Gamma_0)$ is a right-invariant
K\"ahler metric.
\end{corollary}

\begin{remark}
Let $G$ and $\dot{G}$ be Fuchsian groups as in Theorem \ref{Bers}.
Although $T(\dot G)$ and $F(G)$ are submanifolds of $T(\Gamma_0)$
and $F(1)$ respectively, the Bers isomorphism
$\mathfrak{B}:T(\Gamma_0)\rightarrow F(1)$ we study in Section
\ref{Section3} does not induce the Bers isomorphism $\mathfrak{B}_G:
T(\dot{G}) \rightarrow F(G)$. In fact, when $G\neq \{\id\}$, we do
not have an explicit expression of the biholomorphism $h: \U/\dot{G}
\rightarrow (\U/G) \setminus \{\hat{a}\}$ of Theorem \ref{Bers}, and
thus we cannot write down the isomorphism $\mathfrak{B}_G:
T(\dot{G}) \rightarrow F(G)$ and its derivative mapping explicitly
as in Section \ref{section4}. Nevertheless, it will be interesting
to compare the Velling-Kirillov metric on $F(G)$ to the
Takhtajan-Zograf metric on $T(\dot{G})$ under Bers isomorphism
$\mathfrak{B}_G$.
\end{remark}

Bowick and Rajeev \cite{BR1} and Kirillov and Yuriev \cite{KY} have
calculated the curvature tensor of the Velling--Kirillov metric on
the homogenous space $S^1\bk \Diff_+(S^1)$, a natural submanifold of
$S^1\bk\text{Homeo}_{qs}(S^1)\simeq F(1)$. Via the Bers isomorphism
$\mathfrak{B}:T(\Gamma_0)\rightarrow F(1)$, we can use their result
to find the curvature tensor of the Takhtajan-Zograf metric on
$T(\Gamma_0)$ and study the properties of the curvature. This will
in turn give the curvature properties of the Takhtajan-Zograf metric
on the Teichm\"uller spaces of cofinite punctured Riemann surfaces.
We will consider this question in a subsequent paper.

\vspace{1cm} \noindent \textbf{Acknowledgments.} I would like to
thank Leon Takhtajan for the valuable discussions when this work was
done. This work was partially supported by MMU internal funding
PR/2006/0590.

\end{document}